%% file: Gromova_Tur_Barsuk.tex
\documentclass[graybox]{svmult}

\usepackage{mathptmx}       % selects Times Roman as basic font
\usepackage{helvet}         % selects Helvetica as sans-serif font
\usepackage{courier}        % selects Courier as typewriter font
\usepackage{type1cm}        % activate if the above 3 fonts are
\usepackage{makeidx}         % allows index generation
\usepackage{graphicx}        % standard LaTeX graphics tool
\usepackage{multicol}        % used for the two-column index
\usepackage[bottom]{footmisc}% places footnotes at page bottom
\usepackage{multirow}
\usepackage{amsmath}
\makeindex             % used for the subject index
\begin{document}

\title*{A Pollution Control Problem for the Aluminum Production in
Eastern Siberia: Differential Game Approach\thanks{The reported study was funded by RFBR according to the research project  N 18-00-00727 (18-00-00725)}}
 \titlerunning{Differential Game of  Pollution Control  for the Aluminum Production in
Eastern Siberia}
\author{Ekaterina V. Gromova,
Anna V. Tur and Polina I. Barsuk}
\institute{Ekaterina V. Gromova \at St.Petersburg State University,
7/9, Universitetskaya nab., St. Petersburg  199034, Russia
\email{e.v.gromova@spbu.ru} \and Anna V. Tur \at St. Petersburg State
University, St. Petersburg,
Russia  \email{a.tur@spbu.ru} \and Polina I. Barsuk \at
St.Petersburg State University, Saint-Petersburg, Russia \email{polina.barsuk98@gmail.com} }
\maketitle

\abstract*{In this paper, we apply a dynamic game-theoretic model
and analyze the problem of pollution control in Eastern Siberia
region of Russia. When carrying out the analysis we use real
numerical values of parameters. It is shown that cooperation between
the major pollutants can be beneficial not only for the nature but
also for the respective companies.}

\abstract{In this paper, we apply a dynamic game-theoretic model and
analyze the problem of pollution control in Eastern Siberia region
of Russia. When carrying out the analysis we use real numerical
values of parameters. It is shown that cooperation between the major
pollutants can be beneficial not only for the nature but also for
the respective companies.}

\section{Introduction}
\label{intro} Air pollution is a major environmental problem that
affects everyone  in  the civilized world. Emissions from large
industrial enterprises have a great adverse impact on the
environment and the people's quality of life.

%According to the World
%Health Organization, poor condition of the air has a serious
%toxicological impact on human health. More precisely, exposure to
%high levels of air pollution increases the risk of respiratory
%infections, heart disease, stroke and lung cancer \cite{WHO}.
% Also
%climate change is one of the most important problems today.
A
detailed review of the scientific literature published in 1990--2015
on the topic of climate and environmental changes can be found in
\cite{Climate}. Air pollution is closely linked to climate change.
Therefore, one of the most important issues in ecologic management
concerns the reduction of the pollutant emission into the
atmosphere.

Game theory offers a powerful tool for modeling and analyzing
situations where multiple players pursue different but not
necessarily opposite goals. In particular, it is well suited for
analyzing the ecological management problems in which players
(countries, plants) produce some goods while bearing costs due to
the emitted pollution
\cite{ZacZah,Dockner2, Long,Long2,Maler,Ploeg}.
It should be noted that most results on pollution control turn out to be of more theoretical nature because  it is difficult to
obtain realistic numerical values of the model parameters.

In contrast to the mentioned approach we consider local situations
that can be modeled with more precision. Furthermore, we hope that
the obtained results can be of use when planning local policies
aimed at decreasing pollution load in particular regions. Recently,
there has been a paper devoted to the pollution control problem of the city Bratsk from
the Irkutsk region of the Russian Federation based on data for 2011,
\cite{TG}. Our contribution extends the model presented in the
mentioned paper, moreover, the ecological situation is considered for  the largest alumn enterprises of Eastern Siberia  located in
Krasnoyarsk, Bratsk, and Shelekhov  that the largest plants which produce about
70\% of aluminum in Russia has been built.   In the model we include an  absorbtion which is considered for  different  weather conditions. It is known that  the  ecological situation aggravates in the wintertime
on account of the frequent temperature inversions, weak winds, and
fog  \cite{Akhtimankina,Arguchintseva,Avdeeva}.
The problem of pollution control is formulated with a differential game framework and is considered based on data for 2016  \cite{Pr11,Pr12}.

The paper is structured as follows. In section 2, the description of
the  differential game model is presented. In section 2.1, a
non-cooperative solution is found,  the
 Nash equilibrium is considered as an optimality principal.
Section 2.2 deals with cooperative differential game. The numerical
 example of pollution control for the aluminum production in Eastern Siberia is presented in
the section 3.

\section{A Game-Theoretic Model}
\label{sec:1} Consider a game-theoretic model of  pollution control
based on the models \cite{ZacZah, Grom}.  It is
assumed that on the territory of a given region there are $n$
stationary sources of air pollution involved in the game. Each
player has an industrial production site.  Let the production of
each unit is proportional to its pollution $u_i$. Thus, the strategy
of a player is to choose the amount of pollution emitted to the
atmosphere. We assume that the $n$ sources "contribute" to the same stock of
pollution. Denote the stock of accumulated net emissions by $x(t)$.
The dynamics of the stock is given by the following equation with
initial condition:
\begin{equation}\label{g1.2}\dot{x}(t)=\sum^n_{i=1}u_i(t)-\delta x(t), \ t\in [t_0,
T],  \ x(t_0)=x_0,\end{equation} \noindent where $\delta$ denotes
the environment's self-cleaning capacity. Each player $i$ controls
its emission $u_i\in [0, b_i]$, $b_i>0,$ $i=\overline{1,n}$. The
solution will be considered in the class of open-loop strategies
$u_i(t)$.

The net revenue of player $i$ at time instant $t$ is given by
quadratic functional form: $R_i (u_i(t))=u_i(t)
\left(b_i-\frac{1}{2}u_i(t)\right), \ t\in [t_0,T],$ where $b_i>0$.
Each player $i$ bears pollution costs defined as $d_ix(t)$,  where
$d_i\geq 0$ is a fine for environmental pollution. The revenue of
player $i$ at time instant $t$ $R_i (u_i(t))-d_ix(t)$.  The
objective  of player $i$  is to maximize its payoff
\begin{equation}
\label{1.6}K_i\left(x_0,T-t_0,u_1,u_2,\dots,u_n\right)=\int_{t_0}^T\big(R_i(u_i)-d_ix(s)\big)ds.
\end{equation}

\subsection{Nash Equilibrium}
\label{sec:2}

We choose the Nash equilibrium as the principle of optimality in
non-cooperaive game.  To find the  optimal emissions
$u^{NE}_1,\ldots,u^{NE}_n$ for players $1,\ldots,n$, we apply
Pontrygin's maximum principle. The Hamiltonian for this problem is
as follows:
\begin{equation}H_i(x_0,T-t_0,u,\psi)=u_i(t)\left(b_i-\frac{1}{2}u_i(t)\right)-d_ix(t)+\psi_i\left(\sum^n_{i=1}u_i(t)-\delta x(t)\right).\end{equation}
 From the first-order optimality condition
%\begin{equation}
%\frac{\partial H_i(x_0,T-t_0,u,\psi)}{\partial
%u_i}=b_i-u_i(t)+\psi_i(t)=0,\  i=1,\ldots,n\end{equation} \noindent
we get the following formulas for optimal controls: $
u^{NE}_i=b_i+\psi_i, \quad i=1,\ldots,n$.
 Adjoint variables
$\psi_i(t)$ can be found from differential equations $
\frac{\partial H_i(x_0,T-t_0,u,\psi)}{\partial
x}=-\frac{d\psi_i(t)}{dt},\quad \psi_i(T)=0, \quad i=1,\ldots,n$. Then
\begin{equation}\label{8}
u^{NE}_i(t)=b_i-\frac{d_i}{\delta}+\frac{d_i}{\delta}e^{\delta(t-T)},
\quad i=1,\ldots,n.
\end{equation}
Here we assume that for the environment's self-cleaning capacity
$\delta$ the following inequalities hold: $\delta\geq
\frac{d_i}{b_i},$ $i=1,\ldots,n$. This condition ensures that
$u_i\in [0, b_i]$,
 $i=\overline{1,n}$.

 Let $b_N=\sum\limits_{i=1}^nb_i$, $d_N=\sum\limits_{i=1}^nd_i$.
 Then the optimal trajectory is:
\begin{equation}
x^{NE}(t)=C_1e^{-\delta
t}+\frac{d_N}{2\delta^2}e^{\delta(t-T)}+\frac{b_N}{\delta}-\frac{d_N}{\delta^2},
\end{equation}
\noindent where $C_1=e^{\delta
t_0}(x_0-\frac{b_N}{\delta}+\frac{d_N}{\delta^2}-\frac{d_N}{2\delta^2}e^{\delta(t_0-T)})$.

But in the case when for some $i$: $\delta< \frac{d_i}{b_i},$ it may
happens that optimal control $u^{NE}_i$ for player $i$  leaves the
compact $[0; b_i]$. Let
$\overline{t}_i=T+\frac{1}{\delta}ln(1-\frac{b_i\delta}{d_i}).$ If
$\delta< \frac{d_i}{b_i}$ and $\overline{t}_i>t_0$ then optimal
control for player $i$ has a following form:
\begin{equation}u^{NE}_i(t)=\left\{\begin{array}{lr}
        0, & \mbox{for } t_0\leq t\leq \overline{t}_i;\\
         b_i-\frac{d_i}{\delta}+\frac{d_i}{\delta}e^{\delta(t-T)}, & \mbox{for } \overline{t}_i \leq t\leq T.\\
            \end{array}\right.
\end{equation}
It can be noted that
$T+\frac{1}{\delta}ln(1-\frac{b_i\delta}{d_i})\geq
T-\frac{b_i}{d_i}$ for all $\delta> 0$. It means if $T\leq
t_0+\frac{b_i}{d_i},$ then optimal controls have no switching
points.

\subsection{Cooperative Solution}
Consider now the cooperative case of the game. Assume the
players agreed to cooperate and their  goal is to
achieve the joint optimum. The joint payoff is:
\begin{equation}
\sum\limits_{i=1}^nK_i\left(x_0,T-t_0,u_1,u_2,\dots,u_n\right)=\int_{t_0}^T\big(\sum\limits_{i=1}^nR_i(u_i)-d_Nx(s)\big)ds.
\end{equation}

Similarly to non-cooperative case we apply
Pontrygin's maximum principle and obtain:
%
%The Hamiltonian for this problem is as follows:
%\begin{equation}H(x_0,T-t_0,u,\psi)=\sum\limits_{i=1}^nu_i(t)\left(b_i-\frac{1}{2}u_i(t)\right)-d_Nx(t)+\psi\left(\sum^n_{i=1}u_i(t)-\delta x(t)\right).\end{equation}
% From the first-order optimality condition
%\begin{equation}
%\frac{\partial H(x_0,T-t_0,u,\psi)}{\partial
%u_i}=b_i-u_i(t)+\psi(t)=0,\quad i=1,\ldots,n\end{equation} \noindent
%we get the following formulas for optimal controls: $
%u^{*}_i=b_i+\psi, \ i=1,\ldots,n$.
%
% Adjoint variable $\psi(t)$
%can be found from differential equations
%$$ \frac{\partial H(x_0,T-t_0,u,\psi)}{\partial
%x}=-\frac{d\psi(t)}{dt},\quad \psi(T)=0, \quad i=1,\ldots,n. $$
%
%Then

\begin{equation}\label{15}
u^{*}_i(t)=b_i-\frac{d_N}{\delta}+\frac{d_N}{\delta}e^{\delta(t-T)},
\quad i=1,\ldots,n.
\end{equation}
Here we assume that for the environment's self-cleaning capacity
$\delta$ the following inequalities hold: $\delta\geq
\frac{d_N}{b_i},$ $i=1,\ldots,n$. This condition ensures that
$u^{*}_i\in [0, b_i]$,
 $i=\overline{1,n}$.  Then the optimal cooperative trajectory is:
\begin{equation}
x^{*}(t)=C_2e^{-\delta
t}+\frac{nd_N}{2\delta^2}e^{\delta(t-T)}+\frac{b_N}{\delta}-\frac{nd_N}{\delta^2},
\end{equation}
\noindent where $C_2=e^{\delta
t_0}(x_0-\frac{b_N}{\delta}+\frac{nd_N}{\delta^2}-\frac{nd_N}{2\delta^2}e^{\delta(t_0-T)}).$

Notice that the open loop Nash equilibrium  yields more pollution
than the optimal strategies in the cooperative game:
$\sum\limits_{i=1}^nu^{*}_i(T) =\sum\limits_{i=1}^nu^{NE}_i(T)$, and
for all $t\in[t_0;T)$:
$$ \sum\limits_{i=1}^nu^{*}_i(t)
-\sum\limits_{i=1}^nu^{NE}_i(t)=\frac{d_N(n-1)}{\delta}(e^{\delta(t-T)}-1)<0.$$

Also consider the situation, when for  player $i$ the inequality
$\delta< \frac{d_N}{b_i}$ holds. In this case $u^{*}_i(t)$ becomes
negative when $t<\widetilde{t}$,  where
$\widetilde{t}_i=T+\frac{1}{\delta}ln(1-\frac{b_i\delta}{d_N})$. So,
if $\delta< \frac{d_N}{b_i}$ and $\widetilde{t}_i>t_0$, then
 optimal control for player $i$ has a following form:
\begin{equation}u^{*}_i(t)=\left\{\begin{array}{ll}0,    & \mbox{for } t_0\leq t\leq \widetilde{t}_i;\\
  b_i-\frac{d_N}{\delta}+\frac{d_N}{\delta}e^{\delta(t-T)}, & \mbox{for } \widetilde{t}_i \leq t\leq T. \end{array}\right.\end{equation}
It can be noted that
$T+\frac{1}{\delta}ln(1-\frac{b_i\delta}{d_N})\geq
T-\frac{b_i}{d_N}$ for all $\delta> 0$. It means if $ T\leq
t_0+\frac{b_i}{d_N},$ then optimal cooperative control of player $i$
has no switching points.

\section{A Pollution Control Problem in Eastern Siberia}

Non-ferrous metallurgy is one of the most developed industries in
Eastern Siberia. Large aluminum smelters such as Krasnoyarsk, Bratsk
and Irkutsk Aluminum Plants are located in this region. All of the
above-mentioned factories belong to the United Company RUSAL, which
is one of the world's major producers of aluminium.
%However, the
%cities remain on the list of the most polluted ones  in Russia
%despite Rusal's increasing investment in environmental programs.
%%Repeated adverse weather conditions which contribute to the
%accumulation of air pollution exacerbate the situation. According to
%the Article 19 of the Federal Law "On Atmospheric Air Protection"
%No. 96-FZ, legal entities that emit harmful (polluting) substances
%into the air are obliged to take measures concerning reduction of
%the emissions whether it is received forecasts of adverse weather
%conditions.
 In the model under
consideration, the problem of reducing  emissions from smelters
during adverse weather conditions can be solved by changing the
parameter $\delta$ denoted the environment's self-cleaning capacity.

 We consider the 3-players differential game, where  players are
the specified companies. To calculate the required model parameters $b_i$, $d_i$,  we use the
data about the sources of air pollution for year 2016. Let the
coefficient $b_i\geq 0$ equals to a ratio of operating profit of
company  $(P_i)$ to its amount of air emissions $(V_i)$.
Furthermore, $d_i\geq 0$ determines the amount of fine for air
pollution depending on the  total pollution. To determine the fines,
we used the data about companies payments  for air pollution in the
year 2016. Let $L_i$ be the payment  for air pollution of the
company $i$, then:
\begin{equation}\label{b}b_i= \frac{P_i}{V_i}, \quad d_i=\frac{L_i}{V_1+V_2+V_3}. \end{equation}
Table 1 includes the data corresponding to 2016 on the operating
profit of each company, its air pollution and payments for air
pollution.  The operating profit of Krasnoyarsk Aluminum Smelter
could be found in \cite{PrKr}. \cite{Pr1} gives us the joint
operating profit of Bratsk and Irkutsk Aluminum Smelters, which is
equal to 4210,43 million rubles.  We estimated the profit of each
company in proportion to the volume of aluminum produced by these
companies in 2016. According to \cite{BrYear}, Bratsk Aluminum
Smelter produced 1005500 tons of aluminum and  Irkutsk Aluminum
Smelter -- 415400 tons in 2016. So, the operating profit of the two
plants accounts for 2979,51 million rubles and 1230,92 million
rubles respectively. The payment  for air pollution of Krasnoyarsk
Aluminum Smelter amounted to $L_1=87723,95$ thousands rubles
\cite{LKr} in 2016. According to \cite{Pr4} the payment  for air
pollution of the company Irkutsk Aluminum Smelter accounted for
$L_2=18830$ thousands rubles in the same year. Environmental impact
fee including waste disposal fee of Bratsk Aluminum Smelter is
equals to 65278 in 2016 \cite{Pr4}. According to \cite{Pr6} the
payment  for air pollution of Bratsk Aluminum Smelter is
approximately 90 percent of its total environmental impact fee. So,
we estimated its payment for air pollution at
$L_2=0,9\cdot65278=58780,2$ thousands rubles.  Using formulas
(\ref{b})  we get the respective coefficients of the model $b_i$,
$d_i$ (Table 1).
\begin{table}\label{tab1}
\caption{The operating profits, air pollutions and payments for air
pollution of the companies in 2016. The coefficients of the model}
\begin{tabular}{lllllc}
\hline
\textbf{Company} & {$P_i$ {(mln. rubles)}} & {$V_i$ {(tons)}}& $L_i$ { (ths. rubles) $\ $} & {$b_i$} & {$d_i$}\\
\hline Krasnoyarsk  Aluminum Smelter &  3412,23  & 57800&87723,95 $\
$& 59035,12
 & 525,06 \\
Bratsk Aluminum Smelter& 2979,51& 83578,707& 58780,2&35649,15
  &351,64\\
Irkutsk Aluminum Smelter &1230,92 & 25694,1&18830& 47906,72
 &112,71 \\
\hline
\end{tabular}\end{table}

 Table 2 represents the non-cooperative solutions
obtained for some numerical parameters ($t_0 = 0$, $T= 0,4$). We
consider two cases of meteorological conditions, more precisely,
value $\delta =0,02$ corresponds to adverse weather conditions, for
instance, in winter months and $\delta = 0,2$ to normal weather
conditions. The inequalities $T\leq t_0+\frac{b_i}{d_i},$ $ T\leq
t_0+\frac{b_i}{d_N}$ are satisfied  for the chosen parameter values,
thus, optimal cooperative controls of players have no switching
points (we use (\ref{8}), (\ref{15}) to compute the optimal
strategies).
\begin{table}\label{tab3}
\caption{Nash equilibrium strategies. Payoffs of companies  in Nash
equilibrium}
\begin{tabular}{lll}
\noalign{\smallskip}\hline\noalign{\smallskip}
\textbf{Company} & {$u^{NE}_i$, $\delta = 0,02$} & {$u^{NE}_i$, $\delta = 0,2$}\\
\noalign{\smallskip}\hline\noalign{\smallskip}
 KrAS & $32782,12+26253e^{0,02 t-0,008}$&
$56409,82+2625,3e^{0,2 t-0,08}$
\\
BrAS&$18067,15+17582 e^{0,02 t-0,008}$&$33890,95+1758,2e^{0,2
t-0,08}$
 \\
IrAS& $42271,22+5635,5 e^{0,02 t-0,008}$&$47343,17+563,55e^{0,2
t-0,08}$
\\
\noalign{\smallskip}\hline\noalign{\smallskip}
 & $K_i\left(x_0,T-t_0,u^{NE}\right)$, $\delta = 0,02$ & {$K_i\left(x_0,T-t_0,u^{NE}\right)$, $\delta = 0,2$}  \\
\noalign{\smallskip}\hline\noalign{\smallskip}  KrAS  &
$691063605,8-209,19x_0$&$691203820-201,84x_0$
 \\
BrAS &$250177865,6-140,09x_0$&$250271735-135,18x_0$
   \\
IrAS & $457730701,9-44,9x_0$&$457760774,5-43,33x_0$
 \\
\hline
\end{tabular}
\end{table}\vspace*{-20pt}
Table 3 contains the optimal cooperative strategies.
\begin{table}
\caption{Optimal cooperative strategies}
\begin{tabular}{lll}
\noalign{\smallskip}\hline\noalign{\smallskip}
\textbf{Company} & {$u^{*}_i$, $\delta = 0,02$} &{$u^{*}_i$, $\delta = 0,2$ } \\
\noalign{\smallskip}\hline\noalign{\smallskip}  KrAS &
$9564,62+49470,5e^{0,02 t-0,008}$&$54088,07+4947,05e^{0,2 t-0,08}$
\\
BrAS  &$-13821,35+49470,5e^{0,02 t-0,008}$&$30702,1+4947,05e^{0,2
t-0,08}$
 \\
IrAS& $-1563,78+49470,5 e^{0,02 t-0,008}$&$42959,67+4947,0e^{0,2
t-0,08}$
\\
\hline
\end{tabular}
\end{table}

 If we compare   cooperative and non-cooperative  emissions of
players from Table 2 and 3  it is easy to show  that
 the optimal cooperative emissions are less.
%\begin{figure}[t]
%\includegraphics[width=0.49\textwidth]{Upl1}\hfill \includegraphics[width=0.49\textwidth]{Upl2}
%\includegraphics[width=0.49\textwidth]{Upl3}
%\caption{Comparison of cooperative (red dashed line) and
%non-cooperative (blue solid line) emissions}\label{fig1}
%\end{figure}
Table 4 shows  differences between total air pollution in
cooperative and non-cooperative cases and differences between  the
accumulated emissions. The joint cooperative payoffs and its
differences with sum of payoffs in Nash equilibrium  are also  given
in the Table 4.
\begin{table}
\caption{Differences between total air pollution in cooperative and
non-cooperative case and differences between  the accumulated
emissions. Joint cooperative payoff}
\begin{tabular}{llclc}
\noalign{\smallskip}\hline\noalign{\smallskip}
{\textbf{$\delta$}  } &$\sum\limits_{i=1}^nu^{NE}_i(t){-}\sum\limits_{i=1}^nu^{*}_i(t)$ & $x^{NE}(T){-}x^{*}(T)$ & $ \ \ \sum\limits_{i=1}^nK_i\left(x_0,T-t_0,u^*\right)$&$\sum\limits_{i=1}^nK_i\left(u^*\right){-}\sum\limits_{i=1}^nK_i\left(u^{NE}\right)$\\
\noalign{\smallskip}\hline\noalign{\smallskip}
  $\delta = 0,02 \ $&$98941(1-e^{0,02 t-0,008})$&$157,044$&
$1398986922-394,18x_0$&14748,7\\
$\delta = 0,2$&$9894,1(1-e^{0,2 t-0,08})$&$146,2124$
&$1399250309-380,35x_0$&13979,5 \\
\hline
\end{tabular}
\end{table}

 Consider the Shapley value  as an cooperative solution. To calculate it we use a non-standard
method of construction a characteristic function  proposed in
\cite{GrMar}. According to \cite{GrMar} players from coalition $S$
use (obtained earlier) strategies $u^*_S$ from the optimal profile
$u^*$
 and the players from $N \setminus S$ use (obtained
earlier) strategies $u^{NE}_{N \setminus S}$  from the Nash
equilibrium strategies:
\begin{equation}\label{Vnew}
V^\eta(S, \cdot)= \left\{\begin{array}{ll} \quad 0, &
S=\{\emptyset\},
\\
\underset{i \in S}\sum{K_i(\cdot, u^*_S,u^{NE}_{N\setminus S})},&S
\subset N,
\\
\max\limits_{u_1,\ldots, u_n} \sum\limits_{i=1}^n K_i(\cdot,  u_1,
\ldots, u_n),&S=N.
\end{array}\right.
\end{equation}
%This characteristic function  is technically  easier to construct.
Table 5 contains the characteristic function for our example. The
Shapley values  are presented in Table 6. It is also interesting to
see how much each firm benefits from cooperation as compared to a
non-cooperative case. Table 6 shows this difference. We can observe
that it is profitable to the companies to stick to the cooperative
agreement, however to different extent.\vspace*{-5pt}
\begin{table}
\caption{Characteristic function}
\begin{tabular}{lll}\noalign{\smallskip}\hline\noalign{\smallskip}
&$\delta = 0,02$&$\delta = 0,2$\\
\noalign{\smallskip}\hline\noalign{\smallskip}
$V^\eta(\{1\},x_0, T-t_0) $&$691061320,2-209,19x_0$& $691201653,2-201,84x_0$\\
$V^\eta(\{2\},x_0, T-t_0) $&$250173553-140,09x_0$ &
$250267647,2-135,18x_0$
 \\
$V^\eta(\{3\},x_0, T-t_0) $&$457722552,2-44,9x_0$
&$457753050-43,33x_0$
   \\
$V^\eta(\{1,2\},x_0, T-t_0) $&$941245436,6-349,28x_0$ &
$941479313,4-337,02x_0$
 \\
 $V^\eta(\{1,3\},x_0, T-t_0) $ &$1148794743-254,09x_0$& $1148965007-245,17x_0$\\
 $V^\eta(\{2,3\},x_0, T-t_0) $&$707904167-184,99x_0$& $708028338-178,5x_0$\\
\hline
\end{tabular}
\end{table}\vspace*{-20pt}
\begin{table}
\caption{Shapley value. Difference between the Shapley value and the Nash
equilibrium}
\begin{tabular}{lllcc}
\hline
\textbf{Company} & {$Sh_i(x_0,T-t_0)$} & {$Sh_i(x_0,T-t_0)$} &{$Sh_i-K_i\left(u^{NE}\right)$}&{$Sh_i-K_i\left(u^{NE}\right)$}\\
 &{ $\delta = 0,02$}&{ $\delta = 0,2$}& { $\delta = 0,02$}&{ $\delta = 0,2$}\\\noalign{\smallskip}\hline\noalign{\smallskip}  KrAS &
$691072037,4-209,19x_0$&$691211811,9-201,84x_0$& $8431,6$& 7991,9
\\
BrAS &$250182865,8-140,09x_0$&$250276474,5-135,18x_0$&5000,2&4739,5
 \\
IrAS& $457732018,6-44,9x_0$&$457762022,7-43,33x_0$& $1316,7$& 1248,2
\\
\hline
\end{tabular}
\end{table}

The results show that cooperation is beneficial  for all smelters.
 It should be noted that the higher value of the fine $d_i$,
the more profitable the company is to cooperate. In our example
Krasnoyarsk Aluminum Smelter is most motivated for cooperation.

%Total emissions are reduced when players use cooperative strategies.
%Figure \ref{fig4} presents a detailed picture of how cooperation
%leads to the decrease in the amount of emitted pollutants. The
%dependence of the  difference between  the accumulated emissions on
%$\delta $ is shown. Figure \ref{fig4} illustrates the dependence of
%the difference between cooperative and non-cooperative payoffs on
%$\delta $.

%It can be noted that for small values of the parameter
%$\delta $, when unfavorable weather conditions occur, players are
%more motivated for cooperation. We also observe a greater decrease
%of   accumulated emissions under cooperation during adverse weather
%conditions. This shows that the model, which takes into account a
%self-cleaning ability of the atmosphere, allows more effective
%influence on companies to reduce emissions during adverse weather
%conditions.
%%
%\begin{figure}[t]
%\sidecaption[t]
%\includegraphics[width=1\textwidth]{com}
%\caption{Dependence of the differences
%$\sum\limits_{i=1}^nK_i\left(x_0,T-t_0,u^*\right)-\sum\limits_{i=1}^nK_i\left(x_0,T-t_0,u^{NE}\right)$
%and $x^{NE}(T)-x^*(T)$  on $\delta $.}\label{fig4}
%\end{figure}

%\begin{figure}[b]
%\sidecaption
%\includegraphics[width=0.63\textwidth]{dif}
%\caption{Dependence of the difference
%$\sum\limits_{i=1}^nK_i\left(x_0,T-t_0,u^*\right)-\sum\limits_{i=1}^nK_i\left(x_0,T-t_0,u^{NE}\right)$
%on $\delta $.}\label{fig5}
%\end{figure}

\section*{Conclusion} In this paper, we applied game theory
to analyze the problem of pollution control in Eastern Siberia. In
doing so, we considered the real data  for 2016-2018 years obtained
from statistical and accounting reports.
It can be noted that for small values of the absorbtion parameter
$\delta $, when unfavorable weather conditions occur, players are
more motivated for cooperation. We also observe a greater decrease
of   accumulated emissions under cooperation during adverse weather
conditions. This shows that the model, which takes into account a
self-cleaning ability of the atmosphere, allows more effective
influence on companies to reduce emissions during adverse weather
conditions.
%%
%
%
%The obtained result shows
%that the cooperation among the major pollutants can be beneficial
%both for the environment and the companies, especially for those
%with high levels of pollution.

%
%\begin{acknowledgement}
%The reported study was funded by RFBR according to the research
%project  N 18-00-00727 (18-00-00725)
%\end{acknowledgement}

\input{referenc_Gromova_Tur_Barsuk}

\end{document}

%% file: referenc_Gromova_Tur_Barsuk.tex
%%%%%%%%%%%%%%%%%%%%%%%% referenc.tex %%%%%%%%%%%%%%%%%%%%%%%%%%%%%%